\documentclass[10pt]{article}
\usepackage[active]{srcltx}
\usepackage{amssymb,amsmath,amsfonts,bbm,pifont,upgreek,bbold,accents,ulem} 
\usepackage[colorlinks=true]{hyperref}
\hypersetup{urlcolor=blue, citecolor=red}

\renewcommand{\em}{\sl}

\font\teneufm=eufm10
\font\seveneufm=eufm7
\font\fiveeufm=eufm5
\newfam\eufmfam
\textfont\eufmfam=\teneufm
\scriptfont\eufmfam=\seveneufm
\scriptscriptfont\eufmfam=\fiveeufm
\def\eufm#1{{\fam\eufmfam\relax#1}}

\newcommand\beq[1]{ \begin{equation}\label{#1} }
\newcommand{\eeq}{ \end{equation} }
\newcommand{\beqno}{ \[ }
\newcommand{\eeqno}{ \] }
\newcommand\beqa[1]{ \begin{eqnarray} \label{#1}}
\newcommand{\eeqa}{ \end{eqnarray} }
\newcommand{\beqano}{ \begin{eqnarray*} }
\newcommand{\eeqano}{ \end{eqnarray*} }
\newcommand\arr[1]{\left\{\begin{array}{l}#1\end{array}\right.}
\renewcommand{\theequation}{\arabic{section}.\arabic{equation}}

\newtheorem{theorem}{Theorem}[section]
\newtheorem{definition}{Definition}[section]
\newtheorem{proposition}{Proposition}[section]
\newtheorem{lemma}{Lemma}[section]
\newtheorem{sublemma}{Sublemma}[section]
\newtheorem{remark}{Remark}[section]
\newtheorem{notationalremark}{Notational Remark}[section]
\newtheorem{corollary}{Corollary}[section]
\newtheorem{assumption}{Assumption}[section]
\newtheorem{claim}{Claim}[section]
\newtheorem{tools}{$\negsp\negsp$}[subsection]

\newcommand\thm[1]{ \begin{theorem}\label{#1}}
\newcommand\thmtwo[2]{ \begin{theorem}[#1]\label{#2}}
\newcommand\ethm{ \end{theorem} }
\newcommand\dfn[1]{ \begin{definition}\label{#1} \rm}
\newcommand\dfntwo[2]{ \begin{definition}[#1]\label{#2} \rm}
\newcommand\edfn{ \end{definition} }
\newcommand\pro[1]{ \begin{proposition}\label{#1}}
\newcommand\protwo[2]{ \begin{proposition}[#1]\label{#2}}
\newcommand\epro{ \end{proposition} }
\newcommand\lem[1]{ \begin{lemma}\label{#1}}
\newcommand\lemtwo[2]{ \begin{lemma}[#1]\label{#2}}
\newcommand\elem{ \end{lemma} }
\newcommand\sublem[1]{ \begin{sublemma}\label{#1}}
\newcommand\sublemtwo[2]{ \begin{sublemma}[#1]\label{#2}}
\newcommand\esublem{ \end{sublemma} }
\newcommand\rem[1]{ \begin{remark}\label{#1} \rm}
\newcommand\erem{ \end{remark} }
\newcommand\notrem[1]{ \begin{notationalremark}\label{#1} \rm}
\newcommand\enotrem{ \end{notationalremark} }
\newcommand\cor[1]{ \begin{corollary}\label{#1}}
\newcommand\cortwo[2]{ \begin{corollary}[#1]\label{#2}}
\newcommand\ecor{ \end{corollary} }
\newcommand\asmp[1]{ \begin{assumption}\label{#1}}
\newcommand\asmptwo[2]{ \begin{assumption}[#1]\label{#2}}
\newcommand\easmp{ \end{assumption} }
\newcommand\clm[1]{ \begin{claim}\label{#1}}
\newcommand\eclm{ \end{claim} }

\newcommand\equ[1]{{\rm (\ref{#1})}}

\expandafter\chardef\csname pre amssym.def
at\endcsname=\the\catcode`\@
\catcode`\@=11
\def\undefine#1{\let#1\undefined}
\def\newsymbol#1#2#3#4#5{\let\next@\relax
 \ifnum#2=\@ne\let\next@\msafam@\else
 \ifnum#2=\tw@\let\next@\msbfam@\fi\fi
 \mathchardef#1="#3\next@#4#5}
\def\mathhexbox@#1#2#3{\relax
 \ifmmode\mathpalette{}{\m@th\mathchar"#1#2#3}%
 \else\leavevmode\hbox{$\m@th\mathchar"#1#2#3$}\fi}
\def\hexnumber@#1{\ifcase#1 0\or 1\or 2\or 3\or 4\or 5\or 6\or 7\or
8\or
 9\or A\or B\or C\or D\or E\or F\fi}
\ifcase\@ptsize
 \font\tenmsb=msbm10
 \font\sevenmsb=msbm7
 \font\fivemsb=msbm5
\or
 \font\tenmsb=msbm10 scaled \magstephalf
 \font\sevenmsb=msbm7 scaled \magstephalf
 \font\fivemsb=msbm5  scaled \magstephalf
\or
 \font\tenmsb=msbm10 scaled \magstep1
 \font\sevenmsb=msbm7 scaled \magstep1
 \font\fivemsb=msbm5 scaled \magstep1
\fi
\newfam\msbfam
\textfont\msbfam=\tenmsb
\scriptfont\msbfam=\sevenmsb
\scriptscriptfont\msbfam=\fivemsb
\edef\msbfam@{\hexnumber@\msbfam}
\def\Bbb#1{\fam\msbfam\relax#1}
\def\widehat#1{\setboxz@h{$\m@th#1$}%
 \ifdim\wdz@>\tw@ em\mathaccent"0\msbfam@5B{#1}%
 \else\mathaccent"0362{#1}\fi}
\def\widetilde#1{\setboxz@h{$\m@th#1$}%
 \ifdim\wdz@>\tw@ em\mathaccent"0\msbfam@5D{#1}%
 \else\mathaccent"0365{#1}\fi}

\def\RIfM@{\relax\ifmmode}
\def\nonmatherr@#1{\errmessage{\string#1\space allowed only in math mode}}
\def\Bbb{\RIfM@\expandafter\Bbb@\else
 \expandafter\nonmatherr@\expandafter\Bbb\fi}
\def\Bbb@#1{{\Bbb@@{#1}}}
\def\Bbb@@#1{\fam\msbfam\relax#1}
\def\setboxz@h{\setbox\z@\hbox}
\def\wdz@{\wd\z@}
\catcode`\@=\csname pre amssym.def at\endcsname
\newcommand{\ie}{{\rm i.e.\,}}
\newcommand{\eg}{{\it e.g.\,}}
\newcommand{\etc}{{\it etc.\,}}

\newcommand{\Giu}{{\bigskip\noindent}}
\newcommand{\nl}{{\smallskip\noindent}}
\newcommand{\noi}{{\noindent}}

\newcommand{\qed}{\hskip.5truecm
\vrule width 1.7truemm height 3.5truemm depth 0.truemm
\par\Giu}

\newcommand{\negsp}{\hspace{-.09truecm}}  

\newcommand{\dst}{\displaystyle}

\newcommand{\torus}{ {\Bbb T}   }
\renewcommand{\natural}{ {\Bbb N}   }
\newcommand{\real}{ {\Bbb R}   }
\newcommand{\integer}{ {\Bbb Z}   }

\renewcommand{\a }{ {\alpha}   }
\renewcommand{\b}{ {\beta}   }
\newcommand{\g}{ {\gamma}   }
\newcommand{\G}{ {\Gamma}   }

\renewcommand{\k}{ {\kappa}   }
\renewcommand{\l}{ {\lambda}   }
\renewcommand{\L}{ {\Lambda}   }
\newcommand{\m}{ {\mu}   }
\newcommand{\n}{ {\nu}   }

\newcommand{\p}{ {\pi}   }
\renewcommand{\P}{ {\Pi}   }
\renewcommand{\r}{ {\rho}   }
\newcommand{\s}{ {\sigma}   }

\renewcommand{\t}{ {\tau}   }

\renewcommand{\o}{ {\omega}   }
\renewcommand{\O}{ {\Omega}   }

\renewcommand{\Re}{{\, \rm Re\, }}

\newcommand{\cA}{ {\cal A} }

\newcommand{\cT}{ {\cal T} }
\newcommand{\cR}{ {\cal R} }

\newcommand{\cD}{ {\cal D} }

\newcommand\ppu{{ (1) }}
\newcommand\ppd{{ (2) }}
\newcommand\ppt{{ (3) }}

\newcommand\ppk{{ (k) }}

\newcommand\ppj{{ (j) }}

\newcommand\ppn{{ (n) }}

\newcommand\ppi{{ (i) }}

\newcommand\ppo{{ (0) }}

\newcommand\ul{{\uplambda}}
\newcommand\ux{{\upxi}}
\newcommand\uh{{\upeta}}
\newcommand\up{{\rm p}}
\newcommand\uq{{\rm q}}
\newcommand\uz{{\rm z}}

\newcommand\meas{{\, \rm meas\,}}

\newcommand\id{{\, \rm id \,}}

\newcommand\fk{{ k_1 }}
\newcommand\sk{{ k_2 }}
\newcommand\tk{{ k }}
\newcommand\fg{{\g_1}}
\newcommand\pertnorm{{E}}
\newcommand\KAM{{\hat E}}

\begin{document}

\title{Global Kolmogorov tori in the planetary N--body problem. Announcement of result\footnote{ Research supported by ERC Ideas-Project 306414 ``Hamiltonian PDEs and small divisor problems: a dynamical systems approach'' and STAR Project of Federico II University, Naples.}}

\author{  
Gabriella Pinzari\\
\vspace{-.2truecm}
{\footnotesize Dipartimento di Matematica ed Applicazioni ``R. Caccioppoli''}
\\{\footnotesize Universit\`a di Napoli ``Federico II"}
\vspace{-.2truecm}
\\{\footnotesize Monte Sant'Angelo -- Via Cinthia I-80126 Napoli (Italy)}
\vspace{-.2truecm}
\\{\scriptsize gabriella.pinzari@unina.it}
}\date{June 21, 2014}
\maketitle
\vskip.1in
\noi

\begin{abstract}\footnotesize{
We improve a result in \cite{chierchiaPi11b} by proving the existence of a positive measure set of $(3n-2)$--dimensional quasi--periodic motions in the spacial,  planetary $(1+n)$--body problem away from co--planar, circular motions. We also prove that such quasi--periodic motions reach with continuity corresponding $(2n-1)$--dimensional  ones of the planar problem, once the mutual inclinations go to zero (this is related  to a speculation in \cite{arnold63}).
  The main tool is a  full reduction of the SO(3)--symmetry, which retains symmetry by reflections and highlights a quasi--integrable structure, with a small remainder, independently of eccentricities and inclinations.}
\end{abstract}

{
\nl
{\footnotesize{\bf Keywords:} Quasi--integrable structures for perturbed super--integrable systems. N--body problem. Arnold's Theorem on the stability of planetary motions. Multi--scale KAM Theory.  

\nl
{{\bf MSC2000 numbers:}
34D10,  34C20, 70E55,  70F10, 70F15, 70F07,  37J10, 37J15, 37J25, 37J35, 37J40, 70K45

}}

\maketitle
\tableofcontents
\section{Set Up and Background}
\setcounter{equation}{0}
\label{intro}
\renewcommand{\theequation}{\arabic{equation}}

In \cite{arnold63},  V. I. Arnold, partly solving, but undoubtedly clarifying important mathematical  settings of the more than centennial question (going back to the investigations by Sir Isaac Newton, in the XVII century) on the motions of the planetary system,  asserted his  ``Theorem on the stability of planetary motions'' as follows.

\begin{theorem} [{\cite{arnold63, laskarR95, robutel95, herman09, 
fejoz04, pinzari-th09, chierchiaPi11b}}]\label{APT}
In the  many-body problem there exists a set of initial conditions having a 
positive Lebesgue measure and such that, if the initial positions and velocities belong to this set, the distances of the bodies from each other will remain 
perpetually 
bounded, provided the masses of the planets, eccentricities and inclinations are sufficiently small.
 \end{theorem}

\nl
In this paper, we announce an improvement (Theorem \ref{main2} in the next section)  of  Theorem \ref{APT}. 
To present it, we devote this section to a short  survey of related techniques, referring the reader to the aforementioned literature, or  to the review papers \cite{fejoz13, chierchia13, chierchiaPi14} or, finally, to the introduction of \cite{pinzari13} for more details.

%

\vskip.1in
\noi
Consider $(1+n)$ masses in the configuration space $E^3=\real^3$ interacting through gravity. Let such masses be denoted as $m_0$, $\m m_1$, $\cdots$, $\m m_n$, where $m_0$ is a leading mass (``sun'', of ``order one''), while $\m m_1$, $\cdots$, $\m m_n$ are $n$ smaller masses (``planets'', of ``order $\m$'', with $\m$ a very small number). This problem, a sub-problem (usually referred to as ``planetary'' system) of the more general N-body problem,  emulates the solar system; hence,  the study of it has a relevant physical meaning. It is very natural to regard this system (which is Hamiltonian\footnote{I. e., its motions are described by equations of the form $\dst\arr{\dst\dot y^\ppi_j=-\partial_{x^\ppi_j}{\rm H}_{3+3n}(p,q)\\
\dst x^\ppi_j=\partial_{y^\ppi_j}{\rm H}_{3+3n}(p,q)
 }$ where $(p^\ppi, q^\ppi):=(p^\ppi_1, p^\ppi_2, p^\ppi_3, q^\ppi_1, q^\ppi_2, q^\ppi_3$ are  canonical coordinates of the point-mass $i$,  and ${\rm H}_{3+3n}$ is a suitable $(3+3n)$--degrees of freedom Hamilton function, depending on $(p,q)=(p^\ppo,\cdots,p^\ppn, q^\ppo, \cdots, q^\ppn)$.}) as a small perturbation of the leading dynamical problem consisting into the gravitational interaction of the sun separately with each planet. This corresponds to what follows. After letting the system free of the invariance by translations (\ie, eliminating the motion of the sun), 
one can write the $3n$-degrees of freedom Hamiltonian governing the motions of the planets as
\beqa{hel}{\rm H}_{\rm hel}(y, x)&=&\sum_{i=1}^nh_{\rm 2B}^\ppi(y^\ppi,x^\ppi)+\m f_{\rm hel}(y, x)\nonumber\\
&=&\sum_{i=1}^n
\big(\frac{|y^\ppi|^2}{2{\eufm m}_i}-\frac{{\eufm m}_i{\eufm M}_i}{|x^\ppi|}\big)\nonumber\\
&+&\m\sum_{1\le i<j\le n}\big(\frac{y^\ppi\cdot y^\ppj}{m_0}-\frac{m_im_j}{|x^\ppi-x^\ppj|}\big)\ 
\eeqa
where $x^\ppi=(x^\ppi_1, x^\ppi_2,x^\ppi_3)=q^\ppi-q^\ppo$  denote the ``heliocentric distances'', $y^\ppi=(y^\ppi_1, y^\ppi_2,y^\ppi_3)$ their generalized conjugated momenta and ${\eufm m}_i:=\frac{m_0 m_i}{m_0+\m m_i}$, ${\eufm M}_i:={m_0+\m m_i}$ the ``reduced masses''.

  \nl
In order to exploit the integrability of  the ``two-body terms''
$$\dst h_{\rm 2B}^\ppi:=\frac{|y^\ppi|^2}{2{\eufm m}_i}-\frac{{\eufm m}_i{\eufm M}_i}{|x^\ppi|}$$
a natural approach is to put the system in Delaunay\footnote{Delaunay and (see below)  Poincar\'e coordinates are widely described in the literature.  A definition may be found, \eg, in  \cite{chierchiaPi11c, fejoz13b}. Note that $({\rm H}, {\rm h})\in \real^n\times \torus^n$ are denoted as $(\Theta, \theta)$ in \cite{chierchiaPi11c}. {Delaunay--Poincar\'e coordinates were used by several authors, including Arnold, Nekhorossev, Herman, Laskar, Chenciner, F\'ejoz, Robutel, \etc. We shall see below that, due to the proper degeneracy, there is a certain freedom in choosing canonical coordinates for the planetary system. See the definition of {\sl Kepler map} in \S \ref{Kepler maps}.}} coordinates. This is a system of canonical action--angle variables ($(\L,\G,{\rm H},\ell,g,{\rm h})\in \real^{3n}\times \torus^{3n}$), whose r\^ole is the one of transforming (via the Liouville--Arnold Theorem) $h_{\rm 2B}^\ppi$ into ``Kepler form''   , \ie, a function of actions only.  It is well known that, due to the too many integrals of  $h_{\rm 2B}^\ppi$, this integrated form \beq{Kepler}h^\ppi_{\rm K}=-\frac{{\eufm m}_i^3{\eufm M}_i^2}{2\L_i^2}\eeq exhibits a dramatic loss of degrees of freedom: two actions ($\G_i:=|x^\ppi\times y^\ppi|$ and ${\rm H}_i:=x^\ppi_1 y^\ppi_2-x^\ppi_2 y^\ppi_1$ disappear completely. This circumstance is usually called the ``proper degeneracy''.

\nl
 \nl
 Let us denote as
\beq{Del}{\rm H}_{\rm Del}=h_{\rm K}(\L)+\m f_{\rm Del}(\L,\G,{\rm H}, \ell, g, {\rm h})\eeq
 where \beq{Kep}h_{\rm K}(\L_1,\cdots,\L_n):=-\sum_{1\le i\le n}\frac{{\eufm m}_i^3{\eufm M}_i^2}{2\L_i^2}\eeq
 the system \equ{hel} expressed in Delaunay  coordinates. The purpose is to determine a positive measure set of quasi--periodic motions for this system.

 \vskip.1in
\noi
In 1954 A.N. Kolmogorov \cite{kolmogorov54} discovered a breakthrough property of quasi--integrable dynamical systems: for a regular, slightly perturbed system   $${\rm H}(I,\varphi)={\rm h}(I)+\m {\rm f}(I,\varphi)\qquad (I,\varphi)\in A\times \torus^\n$$
where $A\subset \real^\n$ is open,   a great number of quasi--periodic motions $(I_0, \varphi_0 )\to (I_0,\varphi_0+\partial_I{\rm h}(I_0)t)$ of the unperturbed system ${\rm h}$ may be continued in the  dynamics of the perturbed system,  provided the Hessian  $\partial^2_I{\rm h}(I)$ does not vanish identically in $A$. Due to the proper degeneracy, for the planetary system expressed in Delaunay variables \equ{Del}, taking $I:=(\L,\G,{\rm H})$ and $\varphi:=(\ell,g,{\rm h})$, Kolmogorov's non--degeneracy assumption is clearly violated. Despite of this fact, the perturbing function has good parity properties:  Arnold noticed that such parities help in determining a quasi--integrable structure in {\sl all} the variables for the planetary system, as now we explain. 

\vskip.1in
\noi
Following Poincar\'e, one switches from Delaunay coordinates to a new set of canonical coordinates $(\L,\ul,\uh, \ux, \up, \uq)$. These are not in action--angle form, but are in mixed action--angle (the couples $(\L,\l)$) and rectangular form (the $\uz:=(\uh,\ux,\up,\uq)$). The variables $(\L, \l)$ have roughly the same meaning of the $(\L,\ell)$; the $\uz$ are defined in a neighborhood of $\uz=0\in \real^{4n}$ and the vanishing of $(\uh_i, \ux_i)$ or of $(\up_i,\uq_i)$ corresponds to the vanishing of the $i^{\rm th}$ eccentricity, inclination, respectively. 

\nl
Let us denote as
\beq{Poinc}{\rm H}_{\rm P}=h_{\rm K}(\L)+\m f_{\rm P}(\L,\ul,\uz)\qquad \uz=(\uh,\ux,\up,\uq) \eeq
 the system \equ{hel} expressed in Poincar\'e variables.

\nl
Since the perturbation $f_{\rm hel}$ in \equ{hel} does not change under reflection \beq{refl**}(y^\ppi_1,y^\ppi_2,y^\ppi_3, x^\ppi_1, x^\ppi_2,x^\ppi_3)\to ( {\rm r}_1 y^\ppi_1,{\rm r}_2y^\ppi_2,{\rm r}_3y^\ppi_3, {\rm r}'_1x^\ppi_1, {\rm r}'_2x^\ppi_2, {\rm r}'_3x^\ppi_3)\quad {\rm r}_i,\ {\rm r}'_i=\pm 1\eeq and rotation transformations
\beq{rot inv}(y^\ppi_1,y^\ppi_2,y^\ppi_3, x^\ppi_1, x^\ppi_2,x^\ppi_3)\to \Big( {\rm R}(y^\ppi_1, y^\ppi_2, y^\ppi_3), {\rm R}'(x^\ppi_1, x^\ppi_2,x^\ppi_3)\Big)\quad {\rm R},\ {\rm R}'\in {\rm SO}(3)\eeq
and due to the fact that the transformations (respectively, reflections with respect to the coordinate planes and rotation about the $k$--axis)
\beqa{rotations and reflections}
&&\begin{array}{llllllll}
\cR_1^-:\quad &q'{}^\ppi=\big(-x^\ppi_1,\ x^\ppi_2,\ x^\ppi_3\big)\ ,\quad &p'{}^\ppi=\big(y^\ppi_1,\ -y^\ppi_2,\ -y^\ppi_3\big)\\  \\
\cR_2^-:\quad &q'{}^\ppi=\big(x^\ppi_1,\ -x^\ppi_2,\ x^\ppi_3\big)\ ,\quad &p'{}^\ppi=\big(-y^\ppi_1,\ y^\ppi_2,\ -y^\ppi_3\big)\\ \\
\cR_3^-:\quad &q'{}^\ppi=\big(x^\ppi_1,\ x^\ppi_2,\ -x^\ppi_3\big)\ ,\quad &p'{}^\ppi=\big(y^\ppi_1,\ y^\ppi_2,\ -y^\ppi_3\big)\\ \\
\cR_g:\quad &q'{}^\ppi=\big({\cal R}^\ppt_g(x^\ppi_1,\ x^\ppi_2),\ x^\ppi_3\big)\ ,\quad &p'{}^\ppi=\big({\cal R}^\ppt_g(y^\ppi_1,\ y^\ppi_2),\ y^\ppi_3\big)
\end{array}\eeqa
where
$${\cal R}^\ppt_g:=\left(\begin{array}{lll}
\cos g&-\sin g\\
\sin g&\cos g
\end{array}
\right)\qquad g\in \torus$$
 have a nice expression in Poincar\'e variables, respectively,
\beqa{Rg}&&\begin{array}{llllllll}
\cR_1^-:\quad &(\L_i',\ul_i', \uh_i', \ux_i',\up_i',\uq_i')=(\L_i,-\ul_i,\uh_i,-\ux_i, -\up_i,\uq_i)\\  \\
\cR_2^-:\quad &(\L_i',\ul_i', \uh_i', \ux_i',\up_i',\uq_i')=(\L_i,\p-\ul_i,-\uh_i,\ux_i, \up_i,-\uq_i)\\ \\
\cR_3^-:\quad &(\L_i',\ul_i', \uh_i', \ux_i',\up_i',\uq_i')=(\L_i,\ul_i,\uh_i,\ux_i, -\up_i,-\uq_i)\\ \\
\cR_g:\quad &(\L_i',\ul_i', \uh_i', \ux_i',\up_i',\uq_i')=(\L_i,\ul_i+g, \cR^\ppt_{-g}(\uh_i,\ux_i), \cR^\ppt_{-g}(\up_i,\uq_i))
\end{array}\eeqa
one then sees that the averaged (``secular'') perturbation
$$f_{\rm P}^{\rm av}(\L,\uh,\ux,\up,\uq):=\frac{1}{(2\p)^n}\int_{\torus^n}f_{\rm P}(\L,\ul,\uh,\ux,\up,\uq)d\ul$$
enjoys the following symmetries. If we denote
$$t_j:=\frac{\uh_j-{\rm i}\ux_j}{\sqrt2}\quad t_{j+n}:=\frac{\up_j-{\rm i}\uq_j}{\sqrt2}\quad t^*_j:=\frac{\uh_j+{\rm i}\ux_j}{\sqrt2{\rm i}}\quad t^*_{j+n}:=\frac{\up_j+{\rm i}\uq_j}{\sqrt2{\rm i}}$$
and
$$f_{\rm P}^{\rm av}(\L,t,t^*)=\sum_{a, a^*\in \natural^n}{\eufm F}_{a,a^*}(\L)t^\a {t^*}^{\a^*}$$
the Taylor expansion of $f_{\rm P}^{\rm av}$ in powers of $t$, $t^*$, we then have
\begin{proposition}[D'Alembert rules]
\beqa{D'Alembert}
&&f_{\rm P}^{\rm av}(\L,\uh,\ux,\up,\uq)=\arr{
\dst f_{\rm P}^{\rm av}(\L,\uh,-\ux,-\up,\uq)\\ \\
\dst f_{\rm P}^{\rm av}(\L,-\uh,\ux,\up,-\uq)\\ \\
\dst f_{\rm P}^{\rm av}(\L,\uh,\ux,-\up,-\uq)}\nonumber\\\nonumber\\
&& {\eufm F}_{a,a^*}(\L)\ne 0\ \iff\ |a|_1=|a^*|_1\ ,\eeqa 
where $|a|_1:=\sum_{i=1}^n a_i$.
\end{proposition}
By D'Alembert rules one has that the expansion of $f_{\rm P}^{\rm av}$ around $\uz=0$ contains only even monomials and starts with

\beqano
f_{\rm P}^{\rm av}(\L,\uh,\ux,\up,\uq)&=&C_0(\L)+\sum_{1\le i,j\le n}{\cal Q}^{(h)}_{ij}(\L)(\uh_i\uh_j+\ux_i\ux_j)\nonumber\\
&+&\sum_{1\le i,j\le n}{\cal Q}^{(v)}_{ij}(\L)(\up_i\up_j+\uq_i\uq_j)+{\rm O}(\uz^4)
\eeqano
where  $C_0(\L)$, ${\cal Q}^{(h)}_{ij}(\L)$ and ${\cal Q}^{(v)}_{ij}(\L)$ are suitable coefficients, expressed in terms of Laplace coefficients, computed in \cite{laskarR95, herman09, fejoz04}. This expansion shows that the point  $\uz=(\uh,\ux,\up,\uq)=0$ is an {\sl elliptic equilibrium point} for $f_{\rm P}^{\rm av}(\L,\uh,\ux,\up,\uq)$. A natural question is wether, from here, it is also possible to transform $f_{\rm P}^{\rm av}$ into
\beqano
\breve{\rm H}_{\rm P}(\L,\ul,\uz)&=&h_{\rm K}(\L)+\m\breve f_{\rm P}(\L,\ul,\uz)\eeqano
where $\breve f_{\rm P}^{\rm av}$ is in ``Birkhoff normal form'' (hereafter, BNF) of a suitable order (say, of order three).  This means
\beqa{BirkhoffNF}
\breve f_{\rm P}^{\rm av}&=&C_0(\L)+\sum_{i=1}^n\s_i(\L)w_i
+\sum_{i=1}^n\varsigma_i(\L)w_{i+n}+\sum_{r=2}^3 \sum_{1\le i_1\cdots i_r\le 2n}\t_{i_1\cdots i_r}(\L)
w_{i_1}\cdots w_{i_k}
\nonumber\\
&+&{\rm O}(\uz^7)
\eeqa
where $\s_i(\L)$, $\varsigma_i(\L)$ are the eigenvalues of ${\cal Q}^{(h)}(\L)$, ${\cal Q}^{(v)}(\L)$ and, for $1\le i\le n$, $w_i:=\frac{\uh_i^2+\ux_i^2}{2}$, $w_{i+n}:=\frac{\up_i^2+\uq_i^2}{2}$.
Then  Arnold aims  to solve the problem of the proper degeneracy (and hence to prove  Theorem \ref{APT}) by obtaining Kolmogorov full--dimensional tori  bifurcating from the elliptic  equilibrium $\uz=0$, via the following abstract result.
\begin{theorem}[The Fundamental Theorem, \cite{arnold63}]\label{FT} Let
\beq{prop deg}{\rm H}={\rm h}(I)+\m {\rm f}(I, \varphi, u, v)\quad (I,\varphi, u, v)\in A\times \torus^\n\times B\eeq 
  where $A\subset \real^\n$, $B\subset \real^{2\ell}$ are open, $0\in B$, $(I,\varphi)=(I_1,\cdots,I_\n, \varphi_1,\cdots,\varphi_\n)$, $(u,v)=(u_1,\cdots, u_\ell,v_1, \cdots v_\ell)$
be real--analytic and 
\begin{itemize}
\item[{\rm(i)}] $\dst\det \big(\partial^2_I{\rm h}(I)\big)\not\equiv0$;
\item[{\rm(ii)}] $\dst{\rm f^{\rm av}}:=\frac{1}{(2\p)^m}\int_{\torus^m}f(I,\varphi,u,v)d\varphi=\sum_{r=0}^3\sum_{1\le i_1\cdots i_r\le m}\b_{i_1\cdots i_r}(I) w_{i_1}\cdots w_{i_r}+{\rm O}(u,v)^7$, where $w_i:=\frac{u_i^2+v_i^2}{2}$;
\item[{\rm(iii)}] $\dst \det\big(\b_{ij}(I)\big)\not\equiv0$.
\end{itemize}
Then, for any $\k>0$ one can find a number $\varepsilon_0=\varepsilon_0(\k)$ such that, if $0<\varepsilon<\varepsilon_0$ and $0<\m<\varepsilon^8$, the set ${\rm F}_\varepsilon:=A\times \torus^\n\times B^{2\ell}_\varepsilon(0)$ may be decomposed into a set ${\rm F}^*_\varepsilon$ which is invariant for the motions of ${\rm H}$ and a set ${\rm f}_{\varepsilon}$ the measure of which is smaller than $	\k$. More precisely, ${\rm F}^*_\varepsilon$ foliates into $(\n+\ell)$--dimensional invariant manifolds $\{\cT_\o\}_\o$ close to 
$$I_i=I_i^*(\o)\qquad \varphi_i\in \torus\qquad u_j^2+v_j^2=\varepsilon^2 I_j^*(\o)$$ 
where the motion is analytically  conjugated to the linear flow
$$\theta\to \theta+\o t\qquad \theta\in \torus^{\n+\ell}\ .$$
\end{theorem}
Despite this brilliant strategy, Arnold applied Theorem \ref{FT} to the case of the planar three--body problem only, by explicitly checking assumptions (i)--(iii). For the general case, he was aware of some extra--difficulties, about which he
gave just some vague\footnote{We recall, at this respect, the aforementioned contributions by J. Laskar and P. Robutel for the spatial three--body case and by R. Herman and J. F\'eoz for the general case. Note that, while the strategy followed by Laskar and Robutel is intimately related to \cite{arnold63}, the one by Herman and F\'ejoz uses a different KAM scheme with respect to Theorem \ref{FT} and suitable procedures to bypass certain ``degeneracies'' recalled below.}  indications.

\noi
A first problem is represented by the so--called ``secular degeneracies'' :  the ``first order Birkhoff invariants'' $\s_1$, $\cdots$, $\s_n$, $\varsigma_1$, $\cdots$, $\varsigma_n$ satisfy, identically\footnote{Arnold was aware only of the former relation in \equ{secular resonances}. The latter seems to have been noticed, in its full generality, by M. Herman, in the 90s.}, 
\beq{secular resonances}\varsigma_n\equiv 0\quad\sum_{i=1}^n(\s_i+\varsigma_i)\equiv 0\ .\eeq
Such relations are in contrast with usual non--resonance requirements in order to construct BNF \cite{hoferZ94}. But the problem is only apparent. Indeed, it has been recently understood \cite{maligeRL02, chierchiaPi11b} that, by  the symmetry $\cR_g$ in  \equ{Rg},  only resonances $\sum_{i=1}^n(\s_i(\L) k_i+\varsigma_i(\L) k_{i+n})=0$ with $\sum_{i=1}^{2n}k_i=0$ are really important for the construction of BNF, while  resonances  \equ{secular resonances} do not belong to this class. Moreover, in \cite{fejoz04} it has been proved that they are the only ones to be {\sl identically satisfied}; result next improved in \cite{chierchiaPi11b}, where, by direct computation, it has been seen that they   are the only ones to be satisfied in an {\sl open set}: compare item (v)--(a) of Theorem \ref{main1}.

\vskip.1in
\noi
A much more serious problem  is the following\footnote{Proposition \ref{negative} answers, in particular, a question raised by M. R. Herman, who, in \cite{herman09}, declared not to know if the planetary torsion might vanish identically. { More in general, Proposition \ref{negative} generalizes Laplace resonance in \equ{secular resonances} to any order of BNF.} }

\begin{proposition}[Rotational degeneracy \cite{chierchiaPi11c}]\label{negative}
{ For the system \equ{Poinc}, BNF can be constructed up to any prefixed\footnote{Namely, with $3$ replaced by $p$ and ${\rm O}(\uz^7)$ by ${\rm O}(\uz^{2p+1})$ in \equ{BirkhoffNF}.} 
order $p$
 but all the coefficient $\t_{i_1\cdots i_r}(\L)$ of the generic monomial $w_{i_1}\cdots w_{i_r}$ with some of the $i_k$'s equal to $2n$  vanish identically. }
 \end{proposition} 
for which, in particular, the ``torsion'' matrix (the matrix of the second--order coefficients)
$\t=(\t_{ij})$ has an identically vanishing row and column, hence, 
$$\det\t\equiv 0\ .$$
This violates assumption (iii) of Theorem \ref{FT}.

\nl
 However, such negative result, understood only ``a posteriori'', is just the  counterpart of Theorem \ref{main1} below.

\begin{theorem}[{\cite{pinzari-th09, chierchiaPi11b, chierchiaPi11c}}]\label{main1}
It is possible to  determine a global set of canonical coordinates\footnote{RPS stands for ``Regular'', ``Planetary'' and ``Symplectic''. }
\beq{RPS}RPS=(\L,\l,\eta,\xi,p,q)\eeq
which are related to Poincar\'e coordinates $(\L,\ul,\uh,\ux,\up,\uq)$ by
\beqano
&&\L=\L\ ,\quad \ul=\l+\varphi_1(\L,z)\quad \uh_j+{\rm i}\ux_j=(\uh_j+{\rm i}\xi_j)e^{{\rm i}\varphi_2(\L, z)}+{\rm O}({z}^3)\nonumber\\
&& \up=U(\L)p+{\rm O}({z}^3)\quad \uq=U(\L)q+{\rm O}({z}^3)\qquad (*)
\eeqano
where $U(\L)$ is a $n\times n$ unitary matrix, \ie, verifying $U(\L)U^{\rm t}(\L)=\id$ and $\varphi_1$, $\varphi_2$ are suitable functions defined in a global neighborhood of $z=0$, such that\begin{itemize}
\item[{\rm (i)}] $(p_n,q_n)$ are integrals for $f_{\rm RPS}$.
\item[{\rm (ii)}] D'Alembert rules  \equ{Rg} are preserved, and correspond to the  reflections and the rotation in \equ{rotations and reflections}.
In particular,   denoting as  
$${\rm H}_{\rm RPS}(\L,\l,\bar z)=h_{\rm K}(\L)+\m f_{\rm RPS}(\L,\l,\bar z)$$
the system \equ{hel} expressed in the RPS variables, where $\bar z$ denotes $z$ deprived of $(p_n,q_n)$, 
then 
\item[{\rm (iii)}] The point $\bar z=0\in \real^{2n-1}$, which corresponds to the vanishing of all eccentricities and mutual inclinations, is an elliptic equilibrium point for $\bar z\to f^{\rm av}_{\rm RPS}(\L,\bar z)$.  
\item[{\rm (iv)}] For any fixed $p\in \natural$, $p\ge 2$, it is possible to conjugate ${\rm H}_{\rm RPS}$ to
$$\breve{\rm H}_{\rm RPS}(\L,\breve\l,\breve z)=h_{\rm K}(\L)+\m\breve f_{\rm RPS}(\L,\breve\l,\breve z)$$
where 
\beqano
\breve f_{\rm RPS}^{\rm av}(\L,\breve\l,\breve z)&=&C_0(\L)+\sum_{i=1}^n\s_i(\L) \breve w_i
+\sum_{i=1}^{n-1}\varsigma_i(\L)\breve w_{i+n}
\nonumber\\
&+&\sum_{r=2}^p\sum_{1\le i,j\le 2n-1}\t_{i_1\cdots i_r}(\L) \breve w_{i_1}\cdots \breve w_{i_r}+{\rm O}(\breve z^{2p+1})\ .\eeqano

\item[{\rm (v)}] More precisely, for any $p\in \natural$, $a^\ppu_->0$ if $a^\ppn_+:=\infty$, for any $1\le i\le n-1$, it is possible to choice numbers $a^{(i+1)}_+>a^{(i+1)}_-\gg a^\ppi_+$, such that, if $\cA:=\big\{\L=(\L_1,\cdots,\L_n):\ a^\ppi_-\le a^\ppi(\L_i)\le a^\ppi_+\big\}$, then 
\begin{itemize}
\item[{\rm (a)}] $(\s(\L),\bar\varsigma(\L))\cdot k \ne 0$ for any $\L\in\cA$, $k\in \integer^{2n-1}$, $0<|k|_1\le 2p$, $k\ne (1,\cdots,1)$;
\item[{\rm (b)}] $\det\t(\L)\ne 0$ for any $\L\in\cA$.
\end{itemize}
\end{itemize}
\end{theorem}
Clearly, Theorem \ref{main1} below and Theorem \ref{FT} (with $\n:=n$, $\ell:=2n-1$, $I:=\L$, $\varphi:=\breve\l$, $(u,v):=\breve z$) suddenly imply Theorem \ref{APT}, simply replacing\footnote{Substantially, switching from Poincar\'e to RPS  variables  corresponds to replace the $n$ inclinations of the planets with respect to a prefixed frame $(i,j,k)$,  with $(n-1)$ mutual inclinations among te planets 
plus the negligible inclination of the invariable plane with respect to $k$. Recall that the invariable plane is the plane orthogonal to the total angular momentum C.
} ``inclinations'' with ``mutual inclinations'' in the statement. That (*) and (i) imply Proposition \ref{negative} follows by  a classical unicity argument in BNF, suitably adapted to the properly--degenerate case; see \cite{chierchiaPi11c}. 

\nl
We just mention that the variables \equ{RPS} have been obtained  via a suitable ``Poincar\'e regularization'' of a set of action--angle variables $(\L,\G,\Psi,\ell,\g,\psi)$, which we may call ``planetary'' Deprit variables\footnote{The variables $(\L,\G,\Psi,\ell,\g,\psi)$, in such ``planetary form'',  have been rediscovered by the author during her PhD. Note that, apart for few cases \cite{lidovZ76, ferrerO94} of application to the three--body problem, where they reduce  to {the} variables of Jacobi reduction, Deprit variables seem to have remained un-noticed by most. {See also \cite{chierchiaPi11a} for the proof of the symplecticity of $(\L,\G,\Psi,\ell,\g,\psi)$ found in \cite{pinzari-th09}.}
}, since they are in turn easily related to a set of variables  $({\rm R}, \Phi, \Psi, {\rm r}, \varphi, \psi)$ studied in the 80s by F. Boigey and, in their full generality, by A. Deprit \cite{boigey82, deprit83}. The variables  $(\L,\G,\Psi,\ell,\g,\psi)$ ``unfold''  and extend to any $n\ge 2$ a classical  procedure of reduction of the number of degrees of freedom available only for the $n=2$ case and  known since the XIX century, after Jacobi \cite{jacobi1842, radau1868} (often referred to as ``Jacobi reduction of the nodes''). For the relation between the ``original'' Deprit variables $({\rm R}, \Phi, \Psi, {\rm r}, \varphi, \psi)$ and the planetary version $(\L,\G,\Psi,\l,\g,\psi)$ or the relation between the latter and Jacobi reduction of the nodes, see \cite{pinzari13}.

\section{Result}
Clearly, the  {\sl elliptic equilibrium point} of the secular perturbation at the origin plays a fundamental r\^ole in order to determine a quasi--integrable structure in the problem with respect to {\sl all} of its degrees of freedom. As remarked, such equilibrium is determined by the symmetries of the system (\ie, relations \equ{rotations and reflections}). Once the system is put in a set of coordinates such that SO(3)-- symmetry is completely reduced, hence $\cR_g$ will not play its symmetrizing r\^ole anymore,  the elliptic equilibrium, in general  disappears. On the other hand, reducing completely the number of degrees of freedom has its advantages, since it clarifies the structure of phase space and lets the system free of extra--integrals. It is then natural to ask what is the destine of Kolmogorov tori, in such case.

\nl 
 In order to clarify this and other related questions, let us add some more comments.

\begin{itemize}
\item[\tiny\textbullet] The variables \equ{RPS} realize a {\sl partial reduction} of the SO(3)--invariance: in such variables, the system has $(3n-1)$ degrees of freedom, one over the minimum. As said, this is useful in order to describe with regularity the co-inclined, co--circular configuration and to keep the elliptic equilibrium for $\bar z=0$. On the other hand, the fact of having one more degree of freedom than needed implies that possible $(3n-1)$--dimensional resonant tori  corresponding to rotations in the invariable plane of non--resonant  $(3n-2)$--dimensional tori are missed, with subsequent under--estimate  ($\sim\varepsilon^{4n-2}$ instead of $\sim\varepsilon^{4n-4}$) of the measure of the invariant set ${\rm F}_\varepsilon^*$ mentioned in Theorem \ref{FT}. 

\item[\tiny\textbullet] In \cite{chierchiaPi11b} a construction is shown that allows to switch to a ``full reduction'' to  $(3n-2)$ degrees of freedom. Such  procedure is  a bit involved, but allows, at the end to reduce completely the number of degrees of freedom and, simultaneously,  to deal with one only singularity. It generalizes the analogue singularity of Jacobi variables for $n=2$, for which, the planar configuration is not allowed. Therefore,  one has one has to discard a positive measure set in order to stay away from it.  The measure of the invariant set ${\rm F}^*_\varepsilon$ is therefore estimated as $\sim (\varepsilon^{4n-4}-\varepsilon_0^{4n-4})$ with an arbitrary $0<\varepsilon_0<\varepsilon$.

\item[\tiny\textbullet]
The completely reduced variables that are obtained via the full reduction of the previous item for the $n=2$ case are analogue Jacobi's variables (they are not the same) and lead to the same BNF studied in \cite{robutel95}. Differently from what happens for the above discussed case $n=2$, for $n\ge 3$, the full reduction studied in \cite{chierchiaPi11b} looses (besides the $\cR_g$--symmetry in \equ{D'Alembert}) also reflection symmetries and hence the elliptic equilibrium. Such equilibrium needs to be restored via an Implicit Function Theorem procedure, that is successful in the range of small eccentricities and inclinations.

\item[\tiny\textbullet] From the two previous items one has that, while  a   ``continuity''  (letting the inclinations to zero) between $(3n-1)$--dimensional Lagrangian tori of the partially  reduced problem in space  (whose existence has been discussed in \cite{fejoz04, chierchiaPi11b}) and $(2n)$--dimensional Lagrangian tori of the unreduced planar  problem follows from \cite{chierchiaPi11b}, instead, an analogous continuity between $(3n-2)$--dimensional Lagrangian tori of the fully  reduced problem in space  (again discussed in \cite{fejoz04, chierchiaPi11b}) and $(2n-1)$--dimensional Lagrangian tori of the fully planar  problem (discussed in \cite{chierchiaPi11c}) once inclinations go to zero {\sl is naturally expected but, up to now, remains unproved}. Compare also the arguments in \cite{robutel95, fejoz04} on this issue. As mentioned in the previous section, we recall, at this respect, that a controversial (indeed, erroneous) continuity argument between the planar  Delaunay coordinates and the spacial planetary coordinates obtained via Jacobi reduction of the nodes was argued by Arnold \cite{arnold63} in order to infer  non--degeneracy of BNF of the spacial three--body problem. 

\item[\tiny\textbullet] Recall the definitions of ${\rm F}_\varepsilon$, ${\rm F}_\varepsilon^*$ in Theorem \ref{FT}. In both the cases discussed above (partial and full reduction), the ``density'' of ${\rm F}^*_\varepsilon$ inside of ${\rm F}_\varepsilon$, \ie, the ratio
$${\rm d}:=\frac{\meas{\rm F}^*_\varepsilon}{\meas{\rm F}_\varepsilon}$$
goes to one as $\varepsilon\to 0$. That is, {\sl one has to keep more and more close to the co-inclined, co-circular configuration, in order to encounter more and more tori}.  In \cite{chierchiaPi10} it has been proved that one can take
$${\rm d}=1-\sqrt\varepsilon\ .$$
Note in fact that the perturbative technique which leads to Theorem \ref{FT} (or to its improvement discussed in \cite{chierchiaPi10}) is developed with respect to $\varepsilon$, rather than with respect to the initial parameter $\m$ appearing in \equ{hel}. This circumstance is an intrinsic consequence of the fact that the tori obtained via Theorem \ref{FT} bifurcate from the elliptic equilibrium and that, in general, the Birkhoff series 	\equ{BirkhoffNF} diverges.

\item[\tiny\textbullet] In \cite{arnold63} Arnold realized that, in the case of the planar three--body problem the series \equ{BirkhoffNF}  is instead convergent (in this case $f^{\rm av}_{\rm P}$ is integrable). This allows him to prove
$${\rm d}=1-\chi(\m)$$
where $\chi(\m)\to 0$ as $\m\to 0$. For this particular case, the tori do not bifurcate from the elliptic equilibrium, but a different quasi--integrable structure is exploited in \cite{arnold63} (besides also a different perturbative technique  at the place of Theorem \ref{FT}). In \cite{pinzari13}, a slightly weaker result has been proved for the case of the spacial three--body problem and the planar general problem:
$${\rm d}=1-\chi(\m,\a)$$
where $\a$ denotes the maximum semi axes ration and $\chi(\m,\a)\to 0$ as $(\m,\a)\to 0$. Note that for such cases $f^{\rm av}_{\rm P}$ is not integrable.

\item[\tiny\textbullet]
From the astronomical point of view, the investigation mentioned in the two last items is motivated by the fact that, for example, Asteroids  or trans--Neptunian planets exhibit relatively large inclinations or eccentricities. 
From the theoretical point of view, the question is  to understand wether it is possible to find different quasi--integrable structures in the planetary N--body problem besides the one determined by the elliptic equilibrium. 
\end{itemize}

\nl
We prove the following result.

\begin{theorem}\label{main2} Assume that the semi--major axes of the planets are suitably spaced; let $\a$ denote the maximum of such ratios. If $\a$ is small enough and the mass ratio $\m$ is small with respect to some power of $\a$, one can find a number $\varepsilon_0$ and a positive measure set ${\rm F}^*_{\a,\m}\subset {\rm F}:={\rm F}_{\varepsilon_0}$ of Lagrangian, $(3n-2)$--dimensional, Diophantine tori, the density of which in ${\rm F}$ goes to one as $(\a,\m)\to (0,0)$. 
Letting the maximum of the mutual inclinations going to zero,  such $(3n-2)$--dimensional tori are closer and closer to Lagrangian, $(2n-1)$--dimensional, Diophantine tori of the corresponding planar problem.
\end{theorem}

\nl
In the next sections we provide the main ideas behind the proof  of Theorem \ref{main2}. Note that we shall not enter into the (technical) details of the estimate of the density of ${\rm F}^*_{\a,\m}$, for length reasons. 

\section{Tools and Sketch of Proof}

The proof  of Theorem \ref{main2} relies upon four tools.

\subsection{A symmetric reduction of the SO(3)--symmetry} The first tool is a new set of canonical action--angle coordinates which perform a reduction of the total angular momentum in the $(1+n)$--body problem, and, simultaneously, keep symmetry by reflection and are regular for planar motions. Their definition is as follows.

\nl
Let $a^\ppi\in \real_+$,  $P^\ppi\in \real^3$, with $|P^\ppi|=1$, and $e^\ppi$,  denote, respectively,  the {\sl semi--major axis}, the   {\sl direction of the perihelion} and the eccentricity of the $i^{\rm th}$ instantaneous ellipse ${\eufm E}_i$ through $(x^\ppi, y^\ppi)$; let 
 $\cA^\ppi$, with
 $0\le \cA^\ppi\le \cA^\ppi_{\rm tot}=\p (a^\ppi)^2\sqrt{1-(e^\ppi)^2}$, the  area spanned by $x^\ppi$ on ${\eufm E}_i$ with respect to $P^\ppi$ and ${\rm C}^\ppi=x^\ppi\times y^\ppi$ the $i^{\rm th}$ angular momentum.
Define the following partial sums
\beqa{Si}
{\rm S}^\ppj:=\sum_{k=j}^n{\rm C}^\ppk\qquad 1\le j\le n\eeqa
so that
${\rm S}^\ppu:={\rm C}$ is the {\sl total angular momentum}, while ${\rm S}^\ppn={\rm C}^\ppn$.
Define, finally, the following $n$ couples of ${\rm P}$--nodes, $(\widetilde\n_j,\widetilde{\rm n}_j)_{1\le j\le n}$
\beq{good nodes}\widetilde\n_1:=k^\ppt\times {\rm C}\ ,\quad \widetilde{\rm n}_j:={\rm S}^\ppj\times P^\ppj\ ,\quad \widetilde\n_{j+1}:=P^{(j)}\times {\rm S}^{(j+1)}\ ,\quad \widetilde{\rm n}_n:=P^\ppn\eeq
with $1\le j\le n-1$. Then define the coordinates
\beq{P*}{\rm P}_*=(\L, \chi,\Theta,\ell,\k,\vartheta)\eeq
where
\beqano
\begin{array}{lll}
&\dst\L=(\L_1,\cdots,\L_n)\in \real^n\quad &\ell=(\ell_1,\cdots,\ell_n)\in \torus^n\\\\
&\dst\chi=(\chi_0,\bar\chi)\in \real\times \real^{n-1}& \k=(\k_0,\bar\k)\in \torus\times\torus^{n-1}\\\\
&\dst\Theta=(\Theta_0,\bar\Theta)\in \real\times \real^{n-1}&\vartheta=(\vartheta_0,\bar\vartheta)\in \torus\times\torus^{n-1}
\end{array}
\eeqano
with $\bar\chi=(\chi_1,\cdots,\chi_{n-1})$, $\bar\k=(\k_1,\cdots,\k_{n-1})$, $\bar\Theta=(\Theta_1,\cdots,\Theta_{n-1})$, $\bar\vartheta=(\vartheta_1,\cdots,\vartheta_{n-1})$, via the following formulae.

\beqa{belle*}
\begin{array}{llllrrr}
\dst \Theta_{j-1}=\left\{
\begin{array}{lrrr}
\dst{\rm C}_3:={\rm C}\cdot k^\ppt
\\
\\
\dst
{\rm S}^{(j)}\cdot P^{(j-1)}
\end{array}
\right.& \vartheta_{j-1}=\left\{
\begin{array}{lrrr}
\dst\zeta:=\a_{k^\ppt}(k^\ppu, \widetilde\n_1)\qquad& j=1\\
\\
\dst \a_{P^{(j-1)}}(\widetilde{\rm n}_{j-1}, \widetilde\n_{j})&2\le j\le n
\end{array}
\right.\\ \\
\dst\chi_{j-1}:=\left\{
\begin{array}{lrrr}{\rm G}=|{\rm S}^\ppu|
\\
\\
|{\rm S}^\ppj|
\end{array}
\right.
&
\k_{j-1}:=\left\{
\begin{array}{lrrr}{\eufm g}:=\a_{{\rm S}^\ppu}(\widetilde\n_1, \widetilde{\rm n}_1)\qquad& j=1\\
\\
\a_{{\rm S}^\ppj}(\widetilde\n_j, \widetilde{\rm n}_j)&2\le j\le n
\end{array}
\right.
\\ \\
\L_i:={\eufm M}_i\sqrt{{\eufm m}_i a^\ppi}\qquad & \ell_i:=2\p\frac{\cA^\ppi}{\cA_{\rm tot}^\ppi}:= {\rm mean\ anomaly\ of}\ x^\ppi 
\ {\rm on}\ {\eufm E}_i
\end{array}
\eeqa

\nl
Note that \begin{itemize} 
\item[\tiny\textbullet]  The variables \equ{P*}  are very different from the planetary Deprit variables $(\L,\G,\Psi,\ell,\g,\psi)$ mentioned in the previous section. For example, they do not provide the Jacobi reduction of the nodes when $n=2$. Indeed, the definition of \equ{P*} is based on $2n$ nodes \equ{good nodes}, the nodes between the mutual planes orthogonal to ${\rm S}^\ppj$ and $P^{(j)}$ and $P^{(j)}$ and ${\rm S}^{(j+1)}$. Deprit's reduction is instead based on $n$ nodes, the nodes among the planes orthogonal to the   ${\rm S}^\ppj$'s.  {Let us incidentally mention that, for the three--body case ($n=2$), the variables \equ{belle*} are trickily  related to certain canonical variables introduced in \S 2.2 of \cite{pinzari13}. This  relation  will be explained elsewhere.}

\item[\tiny\textbullet]  While, in the case of the variables $(\L,\G,\Psi,\ell,\g,\psi)$, inclinations among the ${\rm S}^\ppj$'s cannot be let to zero, it is not so for the variables \equ{belle*}, where the planar configuration can be reached  with regularity.  And in fact, in the planar case, the change between planar Delaunay variables $(\L, \G, \ell, g)$ and the planar version $(\L, \chi,\ell, \k)$ of \equ{belle*} reduces to
$$\arr{\L=\L\\
\\
\ell=\ell}\qquad\arr{
\chi_{i-1}=\sum_{j=i}^n\G_n\\ \\
\k_{i-1}=g_i-g_{i-1}
}\quad 1\le i\le n$$ 
with $g_0\equiv 0$. Note incidentally that the variables \equ{belle*} are instead singular in correspondence of  the vanishing of the inclinations about ${ P}^\ppj$ and ${\rm S}^{(j)}$ or ${\rm S}^{(j+1)}$ and ${P^{(j)}}$;  configurations with no physical meaning.

\item[\tiny\textbullet] The variables \equ{P*} have in common with the variables $(\L,\G,\Psi, \ell, \g,\psi)$ and the Delaunay variables $(\L,\G,{\rm H}, \ell, g, {\rm h})$ the fact of being singular for zero eccentricities (since in this case the perihelia are not defined). We however give up any attempt to regularize such vanishing eccentricities. The reason is that the Euclidean lengths of the ${\rm C}^\ppj$'s are not\footnote{Indeed, for $1\le j\le n-1$, $|{\rm C}^\ppj|^2=\chi_{j-1}^2+\chi_{j}^2-2\Theta_{j}^2+2\sqrt{(\chi_{j}^2-\Theta_{j}^2)(\chi_{j-1}^2-\Theta_{j}^2)}\cos{\vartheta_{j}}$.} actions (apart for $\chi_{n-1}=|{\rm C}^\ppn|$) and hence the regularization does not seem\footnote{Recall that $e^\ppj=0$ corresponds to $|{\rm C}^\ppj|=\L_j$.} to be (if existing) easy. Note that, since we are interested to high eccentricities motions, we shall have to stay away from these singularities. 

\item[\tiny\textbullet] 
Another remarkable property of the variables \equ{P*}, besides the one of being regular for zero inclinations  is that they retain the symmetry by reflections, as explained in Proposition \ref{symmetry} below. This does not happen for the variables $(\L,\G, \Psi,\ell, \g,\psi)$. As we shall explain better in the next section, such symmetry property plays a r\^ole in order to highlight a global\footnote{ With a remainder  independent of eccentricities and inclinations; compare Proposition \ref{exponential average}.} quasi--integrable structure of ${\rm H}_{\chi_0}$ in \equ{P*Ham} below and, especially, to have an explicit expression of it. 
\end{itemize}
\begin{proposition}\label{symmetry}
The  action--angle  coordinates \equ{P*}
 are canonical. Moreover, letting ${\rm H}_{\chi_0}$ the system \equ{hel}
in these variables, $(\Theta_0,\vartheta_0,\chi_0)$  are integrals of motion for ${\rm H}_{\chi_0}$, which so takes the form
\beq{P*Ham}{\rm H}_{\chi_0}=h_{\rm K}(\L)+\m f_{\chi_0}(\L,\bar\chi,\bar\Theta,\ell,\bar\k,\bar\vartheta)\ .\eeq
Finally, in such variables,  the reflection\footnote{Note that the reflection in \equ{R2-} is slightly different from $\cR_2^-$ in \equ{rotations and reflections}. This is not important, since indeed in \equ{refl**} the signs ${\rm s}_i$ and ${\rm r}_i$ may be chosen independently. } transformation
\beq{R2-}(y_1^\ppi, y_2^\ppi, y_3^\ppi, x_1^\ppi, x_2^\ppi, x_3^\ppi)\to (y_1^\ppi, -y_2^\ppi, y_3^\ppi, x_1^\ppi, -x_2^\ppi, x_3^\ppi)\eeq
is
$$(\L,\chi,\Theta, \ell, \k,\vartheta)\to (\L,\chi,-\Theta, \ell, \k,2\p\integer^n-\vartheta)\ .$$
Therefore, any of the points
\beq{equilibria}(\Theta,\vartheta)=(0, \p k)\quad k\in \{0,1\}^n\qquad \vartheta\mod 2\p\integer^n\eeq
which represents a\footnote{Depending on the signs of the cosines of the mutual inclinations, there are $2^{n-1}$ planar configurations. The one with all the ${\rm C}^\ppi$ parallel and in the same verse corresponds, in the variables \equ{belle*}, to $(\Theta, \vartheta)=\big((0,\cdots, 0),(\p,\cdots,\p)\big)$.} planar configuration, is an equilibrium point for the function $(\bar\Theta,\bar\vartheta)\to f_{\chi_0}(\L,\bar\chi,\bar\Theta,\ell,\bar\chi,\bar\vartheta)$.
\end{proposition}

\subsection{An integrability property}\label{Kepler maps}

The second tool is  an integrability property of the planetary system. To describe it, we generalize a bit the situation, introducing the concept of {\sl Kepler map}.
\vskip.1in
\noi
{\tiny\textbullet} Given $2n$   positive ``mass parameters'' ${\eufm m}_1$, $\cdots$, ${\eufm m}_n$, ${\eufm M}_1$, $\cdots$, ${\eufm M}_n$, a set ${\eufm X}\subset\real^{5n}$  and a bijection
\beqano
\t:\quad {\eufm X}&\to& \big\{({\eufm E}_1,\cdots,{\eufm E}_n)\in (E^3)^n\ ,\ {\eufm E}_i:\ {\rm ellipse}\big\}\nonumber\\
{\rm X}\in {\eufm X}&\to& \big({\eufm E}_1({\rm X}),\cdots,{\eufm E}_n({\rm X})\big)
\eeqano
which assigns to any ${\rm X}\in{\eufm X}$ a $n$--plet  of ellipses $({\eufm E}_1,\cdots,{\eufm E}_n)$  in the Euclidean space $E^3$ with strictly positive eccentricities and having a common focus ${\rm S}$, 
 we shall say that an {injective} map
\beqano
\phi:\quad ({\rm X}, \ell)\in {\cal D}^{6n}:={\eufm X}\times \torus^n\to (y_\phi({\rm X}, \ell),x_\phi({\rm X}, \ell))\in 
(\real^{3})^n\times (\real^{3})^n
\eeqano
is a {\sl Kepler map} if $\phi$
associates to  $({\rm X}, \ell)\in {\eufm X}\times\torus^n$, with $\ell=(\ell_1,\cdots,\ell_n)$ ({\sl mean anomalies})
an element
$$(y_\phi({\rm X}, \ell),x_\phi({\rm X}, \ell))=(y_\phi^\ppu({\rm X}, \ell_1), \cdots, y_\phi^\ppn({\rm X},\ell_n), x_\phi^\ppu({\rm X}, \ell_1),\cdots, x_\phi^\ppn({\rm X}, \ell_n))$$
in the following way.
Letting, respectively, $P_\phi^\ppi({\rm X})$, $a_\phi^\ppi({\rm X})$, $e_\phi^\ppi({\rm X})$ and  $N_\phi^\ppi({\rm X})$   the   direction from ${\rm S}$ to the perihelion, the semi--major axis, the eccentricity and a prefixed direction of the plane of ${\eufm E}_i({\rm X})$, $x_\phi^\ppi({\rm X}, \ell_i)$ are the coordinates with respect to a prefixed orthonormal frame $(i,j,k)$ centered in ${\rm S}$ of the point of ${\eufm E}_i({\rm X})$ such that $\frac{1}{2}a^\ppi_\phi\sqrt{1-(e^\ppi_\phi)^2}\ell_i$ (mod $\p a^\ppi_\phi\sqrt{1-(e^\ppi_\phi)^2}$) is the area spanned from $P_\phi^\ppi({\rm X})$ to $x_\phi^\ppi({\rm X}, \ell_i)$ relatively to the positive (counterclockwise) orientation determined by   $N_\phi^\ppi({\rm X})$ and \beq{yi}
y_\phi^\ppi({\rm X}, \ell_i)={\eufm m}_i \sqrt{\frac{{\eufm M}_i}{(a^\ppi)^3}}\partial_{\ell_i} x^\ppi_\phi({\rm X}, \ell_i)\ .\eeq

\vskip.1in
\noi 
{\tiny\textbullet} A Kepler map will be called {\sl canonical} if any ${\rm X}\in{\eufm X}$ has the form ${\rm X}=({\rm P}, {\rm Q}, \L)$ where $\L=(\L_1,\cdots,\L_n)=({\eufm m}_1\sqrt{{\eufm M}_1a^\ppu_\phi}, \cdots, {\eufm m}_n\sqrt{{\eufm M}_na^\ppn_\phi})$, ${\rm P}=({\rm P}_1,\cdots,{\rm P}_{2n})$, ${\rm Q}=({\rm Q}_1,\cdots,{\rm Q}_{2n})$ and the map
$$(\L,\ell,{\rm P}, {\rm Q})\to (y,x)=(y^\ppu,\cdots,y^\ppn, x^\ppu, \cdots, x^\ppn)$$ preserves the standard 2-form:
$$
\sum_{i=1}^nd\L_i\wedge d\ell_i+\sum_{i=1}^{2n}d{\rm P}_i\wedge d{\rm Q}_i=\sum_{i=1}^n\sum_{j=1}^3 dy^\ppi_j\wedge dx^\ppi_j\ .$$

\vskip.1in
\noi
Examples of canonical Kepler maps are
\begin{itemize}
\item[{\rm (a)}] The map $\phi_{\rm Del}$ which defines the Delaunay variables $(\L,\G,{\rm H},\ell,g,{\rm h})$;
\item[{\rm (b)}] The map $\phi_{\rm Dep}$ which defines the planetary Deprit variables $(\L,\G,\Psi,\ell,\g,\psi)$;
\item[{\rm (c)}] The map $\phi_{\rm P_*}$ which defines the variables ${\rm P}_*=(\L,\chi,\Theta,\ell,\k,\vartheta)$ in \equ{belle*}.
\end{itemize}
\vskip.1in
\noi
{\tiny\textbullet} The following classical relations then hold for (not necessarily canonical) Kepler maps
\beqa{y}
\frac{1}{2\p}\int_\torus \frac{d\ell_i}{|x_\phi^\ppi|} =\frac{1}{a_\phi^\ppi}\ ,\qquad
\frac{1}{2\p}\int_\torus y_\phi^\ppi d\ell_i=0\ ,\quad \frac{1}{2\p}\int_\torus \frac{x_\phi^\ppi}{|x_\phi^\ppi|^3}d\ell_i=0\ .
\eeqa
\vskip.1in
\noi
{\tiny\textbullet} Given   a canonical Kepler map $\phi$, put ${\rm H}_{\phi}:={\rm H}_{\rm hel}\circ\phi$, where ${\rm H}_{\rm hel}$ is as in \equ{hel}. Then
$$
{\rm H}_{\phi}=h_{\rm K}(\L_1,\cdots,\L_n)+\m\,f_{\phi}({\rm X},\ell_1,\cdots,\ell_n)
$$
where $h_{\rm K}$ is as in \equ{Kep}
and 
$$f_\phi({\rm X},\ell_1,\cdots, \ell_n):=\sum_{1\le i<j\le n}(\frac{y_\phi^\ppi({\rm X},\ell_i)\cdot y_\phi^\ppj({\rm X},\ell_j)}{ m_0}-\frac{ m_i m_j}{|x_\phi^\ppi({\rm X},\ell_i)-x_\phi^\ppj({\rm X},\ell_j)|})$$
is the perturbing function \equ{hel} expressed in the variables $(\L,\ell,{\rm P}, {\rm Q})$. 
Imposing a suitable restriction of the the domain so as to  exclude  {\sl orbit collision},  one has that the  {\sl secular $\phi$--perturbing function}, \ie, the average
$$(f_\phi)_{\rm av}({\rm X}):=\frac{1}{(2\p)^n}\int_{\torus^n}f_\phi({\rm X}, \ell_1,\cdots, \ell_n){d\ell_1\cdots d\ell_n}$$
is well defined.
Due to \equ{yi}, the
 ``indirect'' part of the perturbing function, \ie, the term
$\dst \sum_{1\le i<j\le n}$$ y^\ppi_\phi({\rm X},\ell_i)\cdot y_\phi^\ppj({\rm X},\ell_j)/ m_0$ has zero average and hence 
  $(f_{\phi})_{\rm av}$ is just the  average of the Newtonian (or ``direct'') part:
  $$(f_\phi)_{\rm av}({\rm X})=\sum_{1\le i<j\le n}(f_\phi^{(ij)})_{\rm av}
  $$
 with
 $$(f_\phi^{(ij)})_{\rm av}:=-\frac{ m_i m_j}{(2\p)^2}\int_{\torus^2}\frac{d\ell_id\ell_j}{|x_\phi^\ppi({\rm X}, \ell_i)-x^\ppj_\phi({\rm X}, \ell_j)|}\quad i<j\ .$$

\vskip.1in
\noi
{\tiny\textbullet} If we consider the expansion 
\beqno
(f^{(ij)}_\phi)_{\rm av}=(f^{(ij)}_\phi)_{\rm av}^\ppo+(f^{(ij)}_\phi)_{\rm av}^\ppu+(f^{(ij)}_\phi)_{\rm av}^\ppd+\cdots\eeqno
where
\beqano
(f^{(ij)}_\phi)^\ppk_{\rm av}({\rm X}):=-\frac{ m_i m_j}{(2\p)^2}\frac{1}{k!}\frac{d^k}{d\varepsilon^k}\int_{\torus^2}\frac{d\ell_id\ell_j}{|\varepsilon\,x_\phi^\ppi({\rm X}, \ell_i)-x^\ppj_\phi({\rm X}, \ell_j)|}\Big|_{\varepsilon=0}\eeqano
we have that, in this expansion, the two first terms  depend only on $\L_j$. More precisely, due to
\equ{y},
$$(f^{(ij)}_\phi)^\ppo_{\rm av}=-\frac{ m_i m_j}{a^\ppj}\ ,\quad (f^{(ij)}_\phi)^\ppu_{\rm av}=0\ .$$
Therefore, the term
$(f^{(ij)}_\phi)^\ppd_{\rm av}$ carries the first non--trivial information. In the case of the map $\phi=\phi_{\rm P_*}$, we have

\begin{proposition}\label{integrability}
\begin{itemize}
\item[{\rm (i)}] The functions $(f^{(ij)}_{\phi_{\rm P_*}})^\ppd_{\rm av}$ (more in general, $(f^{(ij)}_{\phi_{\rm P_*}})_{\rm av}$) depend only on
  $\L_i$, $\L_j$, $\Theta_{i}$, $\cdots$, $\Theta_{j\wedge(n-1)}$, 
 $\chi_{i-1}$, $\cdots$, $\chi_{j\wedge(n-1)}$,
  $\k_{i}$, $\cdots$, $\k_{j-1}$, $\vartheta_{i}$, $\cdots$, $\vartheta_{j\wedge(n-1)}$ where $a\wedge b$ denotes the minimum of $a$ and $b$. 
\item[{\rm (ii)}] In particular, for any $1\le i\le n-1$, the nearest--neighbor terms $(f^{(i, i+1)}_{\phi_{\rm P_*}})^\ppd_{\rm av}$ (more, in general, $(f^{(i, i+1)}_{\phi_{\rm P_*}})_{\rm av}$) depend only on   $\L_i$, $\L_j$, $\chi_{i-1}$, $\chi_i$, $\chi_{(i+1)\wedge(n-1)}$, $\Theta_{i}$,  $\Theta_{(i+1)\wedge(n-1)}$, $\k_{i}$, $\vartheta_{i}$,  $\vartheta_{{i+1}\wedge(n-1)}$. 
\item[{\rm (iii)}]  The function  $(f^{(n-1,n)}_{\phi_{\rm P_*}})^\ppd_{\rm av}$ 
depends only on $\L_{n-1}$, $\L_n$, $\chi_{n-2}$, $\chi_{n-1}$, $\Theta_{n-1}$,  $\vartheta_{n-1}$,  while it does not depend on  $\k_{n-1}$. Then it  is integrable. 
\item[{\rm (iv)}]  $(f^{(n-1,n)}_{\phi_{\rm P_*}})^\ppd_{\rm av}$ is integrable in the sense of Arnold--Liouville sense: There exists a suitable  global neighborhood $B^2$ of $0\in \real^2$ (where $0$ corresponds to ${\rm C}^{(\n-1)}\parallel{\rm C}^\ppn$),  a set $A\subset \real^4$ and a real--analytic, canonical change of coordinates
\beqano
&&\phi_1:\ \Big((\L_{n-1},\ \L_n,\ \chi_{n-2},\ \chi_{n-1}),\ (\tilde\ell_{n-1}, \tilde\ell_n, \tilde\k_{n-2}, \tilde\k_{n-1}), \ (p_{n-1},\ q_{n-1})\Big) \nonumber\\
&&\to \Big((\L_{n-1},\ \L_n,\ \chi_{n-2},\ \chi_{n-1}),\ (\ell_{n-1}, \ell_n, \k_{n-2}, \k_{n-1}),\ (\Theta_{n-1},\ \vartheta_{n-1})\Big)
\eeqano
defined on $A\times \torus^4\times B^2$ which transforms $(f^{(n-1,n)}_{\phi_{\rm P_*}})^\ppd_{\rm av}$ into a function $h^{(2n+1)}_{\chi_0}$ depending only on $\L_{n-1}$, $\L_n$, $\chi_{n-2}$, $\chi_{n-1}$, $\frac{p_{n-1}^2+q_{n-1}^2}{2}$.
\end{itemize}
\end{proposition}

\vskip.1in
\noi
Note that 
\begin{itemize}
\item[\tiny\textbullet] The  main point of Proposition \ref{integrability} is  that 
the action $\chi_{n-1}=|{\rm C}^\ppn|$ is an integral for $(f^{(n-1,n)}_{\phi_{\rm P_*}})^\ppd_{\rm av}$. Clearly, this is general: whatever is $\phi$,
 $|{\rm C}^\ppj|$ is an integral for $(f^{(ij)}_\phi)^\ppd_{\rm av}$. This fact has been observed firstly, for the case of the three--body problem, in \cite{harrington69}, using Jacobi reduction of the nodes. In that case Harrington observed that $(f^{(12)}_{\phi_{\rm Jac}})^\ppd_{\rm av}$ depends only on $(\L_1,\L_2,\G_1,\G_2, {\rm G}, \g_1)$ and the integrability is exhibited via the couple $(\G_1,\g_1)$. As we already observed, in such case the planetary Deprit variables and the variables obtained by Jacobi reduction of the nodes are  the same. 

\item[\tiny\textbullet] An important issue that is used in the proof of Theorem \ref{main2} (precisely, in order to check certain non--degeneracy assumptions involved in Theorem \ref{KAM} below) is the {\sl effective integration} of  $(f^{(n-1,n)}_{\phi_{\rm P_*}})^\ppd_{\rm av}$. Clearly, in principle, this could be achieved using any of the sets of variables mentioned in the two previous items: planetary Deprit variables $(\L,\G,\Psi,\ell,\g,\psi)$ or the variables \equ{P*}. However, the integration using planetary Deprit variables carries considerable analytic difficulties and has been performed only qualitatively \cite{lidovZ76, ferrerO94}. Using the variables \equ{P*}, such integration can be achieved by a suitable convergent Birkhoff series, exploiting the equilibrium points in \equ{equilibria}. Compare also Proposition \ref{exponential average} and the comments below.
\end{itemize}

\subsection{Global quasi--integrability of the planetary system} The third  tool is the following

\begin{proposition}\label{exponential average}
There exist natural numbers $m$, $\n_1$, $\cdots$, $\n_m$,  with $\n_1+\cdots+\n_m=3n-2>m$ and a positive real number $s$ such that,  if  the semi--major axes of the planets  are suitably spaced and  the maximum semi--axes ratio $\a$ is sufficiently small, for any positive number $\bar K$   sufficiently small  with respect to some positive power of $\a^{-1}$ and any $\m$ small with respect to some power of $\a$, one can find a number $\r(\a,\bar K)$ which goes to zero as a power law with respect to $\a$ and $\frac{1}{\bar K}$, positive numbers  $\g_1$, $\cdots$, $\g_{m}$ depending only on $\a$ and $\m$, a domain ${\rm D}\subset \real^{3n-2}$ a global neighborhood $B^{2(n-1)}$ of $0\in \real^{2(n-1)}$ and, if $C\subset \real^{2n-1}$ is as in Theorem \ref{KAM} with $\n=3n-2$ and $\ell=n-1$,  a real--analytic and symplectic transformation
$$\big((\hat\L,\hat\chi), (\hat\ell, \hat\k), (\hat p,\hat q)\big)\in C_\r\times \torus_{s}^{2n-1}\times B^{2(n-1)}_{\sqrt{2\r}}\to \big((\L,\bar\chi), (\ell, \bar\k),(\bar\Theta, \bar\theta)\big)$$ 
which conjugates the Hamiltonian in \equ{P*Ham} to \beq{quasi integrable}\hat{\rm H}_{\chi_0}=\hat h_{\chi_0}(\hat\L,\hat\chi,\hat p,\hat q)+\m \hat f_{\chi_0}(\hat\L,\hat\chi,\hat\ell,\hat\k,\hat p,\hat q)\eeq 
where $\hat h_{\chi_0}(\hat\L,\hat\chi,\hat p,\hat q)$ depends on $(\hat p_i, \hat q_i)_{1\le i\le n-1}$ only via $\hat J$ $(\hat p,\hat q)$ $:=$ $(\frac{\hat p_1^2+\hat q_1^2}{2}$, $\cdots$, $\frac{\hat p_{n-1}^2+\hat q_{n-1}^2}{2})$ 
and  letting $\o
$ the gradient of $\hat h_{\chi_0}$ with respect to $(\hat\L,\hat\chi, \hat J)$,  then
${\rm D}\supseteq\o^{-1}(\cD^{\bar K, 3n-2}_{\g_1,\cdots,\g_m,\t})\supset\emptyset$.  
Finally, the following holds.
If $L$, $E$, $\hat\r$ are as in Theorem \ref{KAM}, then one can take $\hat\r=\r$, $L=L_0(\a)/\m$, $E=\m E_0(\a)e^{-\bar K s}$, where $L_0(\a)$, $E_0(\a)$ do not exceed some power of $\a^{-1}$.\end{proposition}

\nl
Here are some comments of the proof of Proposition \ref{exponential average}.

\begin{itemize}
\item[\tiny\textbullet] The function $\hat h_{\chi_0}$ is a sum \beq{hat h}\hat h_{\chi_0}=\sum_{i=1}^{2n-1}\hat h^\ppi_{\chi_0}\eeq
where 
$$\hat h^\ppu_{\chi_0}\ ,\ \cdots\ ,\  \hat h^\ppn_{\chi_0}$$
are close to the respective Keplerian terms
$$h^\ppu_{\rm K}\ ,\ \cdots\ ,\ h^\ppn_{\rm K}$$
in \equ{Kepler}, while 
$$\hat h^{(n+1)}_{\chi_0}\ ,\ \cdots\ ,\  \hat h^{(2n+1)}_{\chi_0}$$
are as follows. $\hat h^{(2n+1)}_{\chi_0}$ is close to the function $\m h^{(2n-1)}_{\chi_0}$, where $ h^{(2n-1)}_{\chi_0}$ is defined in the last item of Proposition \ref{exponential average}. For $n\ge 3$ and $2n-2\ge i\ge n+1$, inductively, $\hat h^{(i)}_{\chi_0}$ is as follows. Consider the ``projection\footnote{By ``projection over normal modes'' of a given function $f(I,\varphi,p,q)=\sum_{(a,b)\in \natural^m\times \natural^m, k\in \integer^n}f_k(I) e^{{\rm k}k\cdot\varphi}\prod_{i=1}^m(\frac{u_i-{\rm i}v_i}{\sqrt2})^{a_i}(\frac{u_i+{\rm i}v_i}{{\rm i}\sqrt2})^{b_i}$ we mean the function $\sum_{a\in \natural^m}f_0(I) \prod_{i=1}^m(\frac{u_i^2+v_i^2}{2{\rm i}})^{a_i}$.} over normal modes'' of $(f^{(i-n, i-n+1)}_{\phi_{\rm P_*}})^\ppd_{\rm av}\circ\phi_1\circ\cdots\circ\phi_{2n-1-i}$ with respect to the variables $(p_j, q_j)$ with $j\ge i-n+1$ and $(\chi_i,\tilde\k_i)$ with $i\ge i-n$. This is a function of $$\L_{i-n}, \cdots, \L_n, \chi_{i-n-1}, \cdots, \chi_{n-1}, \Theta_{i-n}, \vartheta_{i-n}, \frac{p_{i-n+1}^2+q_{i-n+1}^2}{2},\cdots ,\frac{p_{n-1}^2+q_{n-1}^2}{2}$$ and  is integrable in the sense of Liouville--Arnold: there exists  $\phi_{2n-i}$ which lets this projection into a function $h^{(i)}_{\chi_0}$ of $$\L_{i-n}, \cdots, \L_n, \chi_{i-n-1}, \cdots, \chi_{n-1},  \frac{p_{i-n}^2+q_{i-n}^2}{2},\cdots, \frac{p_{n-1}^2+q_{n-1}^2}{2}\ .$$
Then $\hat h^{(i)}_{\chi_0}$ is close to $\m h^{(i)}_{\chi_0}$.

\item[\tiny\textbullet] The exponential decay of $E$ with respect to $\bar K$ follows from a suitable averaging technique  derived from \cite{poschel93}, carefully adapted to our case.

\item[\tiny\textbullet] The functions in \equ{hat h} are of different strengths, with respect to the mass parameter $\m$ and the semi--mjor axes  ratios $\a_i:=\frac{a^\ppi}{a^{(i+1)}}$. The first $n$ ones, which are, as said, close to be Keplerian, are of order
$$\sim\ \frac{1}{a^\ppu}\ ,\ \cdots\ ,\ \frac{1}{a^\ppn}\ .$$
The remaining $(n-1)$ ones are much smaller
$$\sim\ \m\frac{(a^\ppu)^2}{(a^\ppd)^3}\ ,\ \cdots\ ,\ \m\frac{(a^{(n-1)})^2}{(a^\ppn)^3}$$
(they have the strength of $\m (f_{\phi_{\rm P_*}}^{(1,2)})^\ppd_{\rm av}$, ..., $\m (f_{\phi_{\rm P_*}}^{(n-1,n)})^\ppd_{\rm av}$, which are so). Therefore in order to apply a KAM scheme to the Hamiltonian \equ{quasi integrable}, we need a formulation suitably adapted to this case. This is given in the following section.
\end{itemize}

\subsection{Multi--scale KAM theory}
The fourth tool is a multi--scale KAM Theorem. To quote it, let us fix the following notations.

\vskip.1in
\noi
Given $m$, $\n_1$, $\cdots$, $\n_m\in \natural$, $\n:=\n_1+\cdots+\n_m$, let us decompose
$$\integer^\n\setminus\{0\}=\bigcup_{i=1}^{m}{\eufm L}_i\setminus{\eufm L}_{i-1}$$
where 
$$\integer^\n=:{\eufm L}_0\supset{\eufm L}_1\supset{\eufm L}_2\supset \cdots \supset {\eufm L}_m=\{0\}$$
is a decreasing sequence of sub--lattices  defined by
\beqno
{\eufm L}_i:=\big\{k=(k_1,\cdots, k_{m})\in \integer^\n=\integer^{\n_1}\times\cdots\times\integer^{\n_m} : \quad k_1=\cdots=k_{i}=0\big\}\ .\eeqno
Next, given $\g$, $\g_1$, $\cdots$, $\g_m$, $\t\in \real_+$, define the ``multi--scale Diophantine'' number sets
\beqano
&&{\cal D}^{\n, K, i}_{\g;\t}:=\Big\{\o\in \real^\n:\ |\o\cdot k|\ge \frac{\g}{|k|^\t}\quad \forall k\in {\eufm L}_{i-1}\setminus{\eufm L}_i,\ |k|_1\le K\Big\}\nonumber\\
&&{\cal D}^{\n, K}_{\g_1\cdots\g_m;\t}:=\bigcap_{i=1}^m
{\cal D}^{\n, K, i}_{\g_i;\t}\qquad {\cal D}^{\n}_{\g_1\cdots\g_m;\t}:=\bigcap_{K\in \natural}
{\cal D}^{\n, K}_{\g_1\cdots\g_m;\t}
\ .
\eeqano 
Explicitly, a number $\o=(\o_1,\cdots,\o_m)\in\real^\n=\real^{\n_1}\times\cdots\times \real^{\n_m}$ belongs to ${\cal D}^{\n}_{\g_1\cdots\g_m;\t}$ if, for any $k=(k_1,\cdots,k_m)\in 
\integer^{\n_1}\times\cdots\times \integer^{\n_m}\setminus\{0\}$,
 \beqno
|\sum_{j=1}^m\o_j\cdot k_j|\geq
\left\{ \begin{array}{l}
\dst\frac{\fg}{|\tk|^{\t}}\quad {\rm if}\quad \fk\neq0\ ;\\ \ \\
\dst \frac{\g_2}{|k|^{\t}}\quad {\rm if}\quad \fk= 0\ ,\quad \sk\neq 0\ ;\\ \ \\
\cdots\\ \ \\
\dst \frac{\g_m}{|k_m|^{\t}}\quad {\rm if}\quad k_1=\cdots=k_{m-1}= 0,\ \cdots,\ k_{m}\neq 0\ .
\end{array}\right.
\eeqno
 Note that the choice
 $m=1$ gives the usual Diophantine set $\cD^\n_{\g_1,\t}$. The  $m=2$-case,  with $\g_1={\rm O}(1)$ and $\g_2={\rm O}(\m)$ has been considered in \cite{arnold63} (and \cite{chierchiaPi10}) for the proof of Theorem \ref{FT}.

\begin{theorem}[Multi--scale KAM Theorem]\label{KAM}
Let $m$, $\ell$, $\n_1$, $\cdots$, $\n_m\in \natural$, with $\n:=\n_1+\cdots+\n_m\ge \ell$, $\t_*>\n$, $\g_1\ge \cdots\ge\g_m>0$, $0<4s\leq \bar{s}<1$, $\r>0$, ${D}\subset\real^{\n-\ell}\times \real^\ell
$,  $A:={D}_\r$, $B^{2\ell}$ a  neighborhood (with possibly different radii) of $0\in \real^{2\ell}$ such that, if $\bar I(u,v):=(\frac{u_1^2+v_1^2}{2},\cdots, \frac{u_\ell^2+v_\ell^2}{2})$, then $\P_{\real^\ell}{\rm D}=\bar I(B^{2\ell})$, $C:=\P_{\real^{\n-\ell}}D$ and let
\beqano
{\rm H}
(I,\varphi, u, v)={\rm h}(I, u,v)+{\rm f}(I,\varphi, u,v)
\eeqano
be real--analytic on $C_\r\times \torus_{\bar{\eufm s}+s}^{\n-\ell}\times B^{2\ell}_{\sqrt{2\r}}$, where ${\rm h}$ depends on $(u,v)$ only via $\bar I(u,v)$.
Assume that 
 $\o_0:=\partial_{(I,\bar I)} {\rm h}$
is a diffeomorphism of $A$ with non singular Hessian matrix $U:=\partial^2_{(I,\bar I)}{\rm h}$  and let $U_k$ 
denote the $ (\n_k+\cdots+\n_m)\times \n$ 
submatrix of $U$, \ie, the matrix with entries $(U_k)_{ij}=U_{ij}$, for $\n_{1}+\cdots+\n_{k-1}+1\leq i\leq \n$, $1\leq j\leq \n$
, where  $2\le k\le m$.
Let
\beqano
&&  {\rm M}\geq\sup_{A}\|U\|\ ,\quad  {\rm M}_k\geq\sup_{A}\|U_k\|\ ,\quad  \bar {\rm M}
\geq\sup_{A}\|U^{-1}\|\ ,\quad
\pertnorm\geq\|{\rm f}\|_{\r,\bar{\eufm s}+s}\nonumber\\
&&\bar {\rm M}_k\geq \sup_{A}\|T_k\|\quad {\rm if}\quad \dst U^{-1}=\left(\begin{array}{lrr}
T_1\\
\vdots\\
T_m
\end{array}
\right)\qquad 1\le k\le m\ .\eeqano
Define
\beqano
&&  \dst K:=\frac{6}{s}\ \log_+{\left(\frac{\pertnorm {\rm M}_1^2\,L}{\gamma_1^2}\right)^{-1}}\quad {\rm where}\quad \log_+ a :=\max\{1,\log{a}\}\\
&&  \dst \hat\r_k:=\frac{\g_k}{3{\rm M}_kK^{\t_*+1}}\ ,\quad  \hat\r:=\min\left\{\hat\r_1\ ,\ \cdots,\ \hat\r_m,\ 
 \r\right\}\\
 \\
&&  \dst L:=\max\ \Big\{\bar {\rm M}
, \  {\rm M}_1^{-1},\ \cdots,\ {\rm M}_m^{-1}\Big\}
\\
&& \hat E:=\frac{E L}{\hat\r^2}
\eeqano
Then one can find two numbers $\hat c_\n>c_\n$ depending only on $\n$ such that, if
the perturbation ${\rm f}$ so small that the following ``KAM condition'' holds
\beqno
\hat c_\n\KAM<1\ ,
\eeqno
then, for any $\o\in\O_*:=\o_0({ D})\cap\cD^\n_{\g_1,\cdots,\g_m,\t_*}$, one can find a unique real--analytic embedding
\beqano
\phi_\o:\quad 
\vartheta=(\hat\vartheta,\bar\vartheta)\in\torus^{\n
}
&\to&(\hat v(\vartheta;\o),\hat\vartheta+\hat u(\vartheta;\o), \cR_{\bar\vartheta+\bar u(\vartheta;\o)}w_1,\ \cdots,\  \cR_{\bar\vartheta+\bar u(\vartheta;\o)}w_\ell)\nonumber\\
&&\in \Re C_r\times \torus^{\n-\ell}\times \Re B^{2\ell}_{\sqrt{2r}}
\eeqano
where $r:=
c_\n
\KAM \hat\r$  such that ${\rm T}_\o:=\phi_{\o}(\torus^\n)$ is a real--analytic $\n$--dimensional ${\rm H}$--invariant torus, on which the  ${\rm H}$--flow is analytically conjugated to $\vartheta\to \vartheta+\o\,t$. 
\end{theorem}

\nl
Theorem \ref{KAM} is essentially  Proposition 3 of \cite{chierchiaPi10} suitably adapted to our case. Applying Theorem \ref{KAM} to the Hamiltonian \equ{quasi integrable} (with $I:=(\hat\L,\hat\chi)$, $\varphi:=(\hat\ell,\hat\k)$, $(u,v):=(\hat p,\hat q)$, $\n=3n-2$, $\ell=n-1$, $m$, $\n_1$, $\cdots$, $\n_m$ as in Proposition \ref{exponential average}) gives the proof of Theorem \ref{main2}. More details will be published elsewhere. \qed
{\small {\bf Acknowledgments.} The author is indebted to J. F\'ejoz, from whom she learned, in January of 2013, during a conference in Banff,  about the integrability of the first non--trivial term of the secular perturbing function of the three--body problem, known since \cite{harrington69}. This was already used in \cite{pinzari13} (paper firstly submitted in July of 2013) to prove the independence of the invariant set on eccentricities for planetary problem,  in the   three--body case, result  also proved in \cite{pinzari13} for the  general planar one (for which case \cite{harrington69} does not apply directly). She is also grateful to him  for appointing her (in June of 2014, after this note was completed) that, in more than one occasion since the end of 2013 (\eg, CELMEC VI, San Martino al Cimino, Viterbo, Sept. 2013), he talked about the independence of such invariant set for the general problem on planetary masses.
}

{\footnotesize
\bibliographystyle{plain}

\begin{thebibliography}{10}

\bibitem{arnold63}
V.I. Arnold.
\newblock {S}mall denominators and problems of stability of motion in classical
  and celestial mechanics.
\newblock {\em Russian Math. Surveys}, 18(6):85--191, 1963.

\bibitem{boigey82}
F.~Boigey.
\newblock {\'E}limination des n\oe uds dans le probl{\`e}me newtonien des
  quatre corps.
\newblock {\em Celestial Mech.}, 27(4):399--414, 1982.

\bibitem{chierchia13}
L.~Chierchia.
\newblock {T}he {P}lanetary {N}--{B}ody {P}roblem.
\newblock {\em UNESCO Encyclopedia of Life Support Systems}, 6.119.55, 2012.

\bibitem{chierchiaPi10}
L.~Chierchia and G.~Pinzari.
\newblock Properly--degenerate {KAM} theory (following {V}.{I}. {A}rnold).
\newblock {\em Discrete Contin. Dyn. Syst. Ser. S}, 3(4):545--578, 2010.

\bibitem{chierchiaPi11a}
L.~Chierchia and G.~Pinzari.
\newblock {D}eprit's reduction of the nodes revised.
\newblock {\em Celestial Mech.}, 109(3):285--301, 2011.

\bibitem{chierchiaPi14}
L.~Chierchia and G.~Pinzari.
\newblock {M}etric stability of the planetary n--body problem.
\newblock {\em {P}roceedings of the International Congress of Mathematicians},
  2014.

\bibitem{chierchiaPi11c}
Luigi Chierchia and Gabriella Pinzari.
\newblock Planetary {B}irkhoff normal forms.
\newblock {\em J. Mod. Dyn.}, 5(4):623--664, 2011.

\bibitem{chierchiaPi11b}
Luigi Chierchia and Gabriella Pinzari.
\newblock The planetary {$N$}-body problem: symplectic foliation, reductions
  and invariant tori.
\newblock {\em Invent. Math.}, 186(1):1--77, 2011.

\bibitem{deprit83}
A.~Deprit.
\newblock Elimination of the nodes in problems of {$n$} bodies.
\newblock {\em Celestial Mech.}, 30(2):181--195, 1983.

\bibitem{fejoz04}
J.~F{\'e}joz.
\newblock D{\'e}monstration du `th{\'e}or{\`e}me d'{A}rnold' sur la
  stabilit{\'e} du syst{\`e}me plan{\'e}taire (d'apr{\`e}s {H}erman).
\newblock {\em Ergodic Theory Dynam. Systems}, 24(5):1521--1582, 2004.

\bibitem{fejoz13}
J.~F{\'e}joz.
\newblock On "{A}rnold's theorem" in celestial mechanics --a summary with an
  appendix on the poincar{\'e} coordinates.
\newblock {\em Discrete and Continuous Dynamical Systems}, 33:3555--3565, 2013.

\bibitem{fejoz13b}
Jacques F{\'e}joz.
\newblock On action-angle coordinates and the {P}oincar\'e coordinates.
\newblock {\em Regul. Chaotic Dyn.}, 18(6):703--718, 2013.

\bibitem{ferrerO94}
Sebasti{\'a}n Ferrer and Carlos Os{\'a}car.
\newblock Harrington's {H}amiltonian in the stellar problem of three bodies:
  reductions, relative equilibria and bifurcations.
\newblock {\em Celestial Mech. Dynam. Astronom.}, 58(3):245--275, 1994.

\bibitem{harrington69}
Robert~S. Harrington.
\newblock The stellar three-body problem.
\newblock {\em Celestial Mech. and Dyn. Astrronom}, 1(2):200--209, 1969.

\bibitem{herman09}
M.~R. Herman.
\newblock Torsion du probl{\`e}me plan{\'e}taire, edited by J. F{\'e}joz in
  2009.
\newblock Available in the electronic `Archives Michel Herman' at
  {\small\url{http://www.college-de-france.fr/default/EN/all/equ_dif/archives_michel_herman.htm}}.

\bibitem{hoferZ94}
H.~{H}ofer{,} E.~Zehnder.
\newblock {\em {S}ymplectic {I}nvariants and {H}amiltonian {D}ynamics}.
\newblock Birkh{{\"a}}user {V}erlag, Basel, 1994.

\bibitem{jacobi1842}
C.~G.~J. Jacobi.
\newblock Sur l'{\'e}limination des noeuds dans le probl{\`e}me des trois
  corps.
\newblock {\em Astronomische Nachrichten}, Bd XX:81--102, 1842.

\bibitem{kolmogorov54}
A.N. Kolmogorov.
\newblock On the {C}onservation of {C}onditionally {P}eriodic {M}otions under
  {S}mall {P}erturbation of the {H}amiltonian.
\newblock {\em Dokl. Akad. Nauk SSR}, 98:527--530, 1954.

\bibitem{laskarR95}
J.~Laskar and P.~Robutel.
\newblock Stability of the planetary three-body problem. {I}. {E}xpansion of
  the planetary {H}amiltonian.
\newblock {\em Celestial Mech. Dynam. Astronom.}, 62(3):193--217, 1995.

\bibitem{lidovZ76}
M.~L. Lidov and S.~L. Ziglin.
\newblock Non-restricted double-averaged three body problem in {H}ill's case.
\newblock {\em Celestial Mech.}, 13(4):471--489, 1976.

\bibitem{maligeRL02}
F.~Malige, P.~Robutel, and J.~Laskar.
\newblock Partial reduction in the {$n$}-body planetary problem using the
  angular momentum integral.
\newblock {\em Celestial Mech. Dynam. Astronom.}, 84(3):283--316, 2002.

\bibitem{pinzari-th09}
G.~Pinzari.
\newblock {\em {O}n the {K}olmogorov set for many--body problems}.
\newblock PhD thesis, {U}niversit{\`a} {R}oma {T}re, April 2009.

\bibitem{pinzari13}
Gabriella Pinzari.
\newblock Aspects of the planetary {B}irkhoff normal form.
\newblock {\em Regul. Chaotic Dyn.}, 18(6):860--906, 2013.

\bibitem{poschel93}
J.~P{\"o}schel.
\newblock Nekhoroshev estimates for quasi-convex {H}amiltonian systems.
\newblock {\em Math. Z.}, 213(2):187--216, 1993.

\bibitem{radau1868}
R.~Radau.
\newblock Sur une transformation des {\'e}quations diff{\'e}rentielles de la
  dynamique.
\newblock {\em Ann. Sci. Ec. Norm. Sup.}, 5:311--375, 1868.

\bibitem{robutel95}
P.~Robutel.
\newblock Stability of the planetary three-body problem. {II}. {KAM} theory and
  existence of quasiperiodic motions.
\newblock {\em Celestial Mech. Dynam. Astronom.}, 62(3):219--261, 1995.

\end{thebibliography}

\def\cprime{$'$} \def\cprime{$'$}

\addcontentsline{toc}{section}{References}
}
\end{document}